\documentclass{amsart}
\usepackage{amsmath, amssymb}
\newtheorem{Theorem}{Theorem}[section]
\newtheorem{Cor}[Theorem]{Corollary}

\newtheorem{Proposition}[Theorem]{Proposition}
\theoremstyle{definition}
\newtheorem{Definition}[Theorem]{Definition}
\theoremstyle{remark}
\newtheorem{rem}[Theorem]{Remark}
\numberwithin{equation}{section}

\newcommand{\R}{\mathbb R}
\newcommand{\N}{\mathbb N}
\newcommand{\Z}{\mathbb Z}

\newcommand{\h}{\mathcal{H}}
\newcommand{\C}{\mathcal{C}}
\newcommand{\F}{\mathcal{F}}
\newcommand{\p}{\mathcal{P}}
\newcommand{\A}{\mathcal{A}}

\newcommand{\rel}{\mathcal{R}}

\newcommand{\mch}{{\mathcal H}}

\begin{document}

\title[Series over a groupoid and cobordism]{On algebras  and groups of formal series over a groupoid and application to some spaces of cobordism}
\author{Jean-Pierre Magnot}

\address{LAREMA, Universit\'e d'Angers, 2 Bd Lavoisier, 49045 Angers cedex 1, ´
	France and Lyc\'ee Jeanne d'Arc \\ Avenue de Grande Bretagne \\ F-63000 Clermont-Ferrand}
\email{jean-pierr.magnot@ac-clermont.fr}


\begin{abstract}
We develop here a concept of deformed algebras and their related groups through two examples. Deformed algebras are obtained from a fixed algebra by deformation along a family of indexes, through formal series. We show how the example of deformed algebra used in \cite{Ma2013} is only an example among others, and how they often give rise to regular Fr\"olicher Lie groups, and sometimes to Fr\'echet Lie groups, that are regular. The first example, indexed by $\N,$ is obtained by deformations of the gourp of bounded classical pseudo-differential operators $Cl^{0,*}$ by algebras of (maybe unbounded) classical pseudo-differential operators. In the second one, the set of indexes is a $\N-$graded groupo\"id, which is made of manifolds with boundary that are understood as morphisms of cobordisms. Here again, we get regular Lie groups, and we show how this setting applies to a class of examples that are derived of so-called stochastic cosurfaces.  
\end{abstract}
\maketitle

\tableofcontents
\section*{Introduction}
In \cite{Ma2013}, an algebra and a group of formal series of operators is described in order to rewrite the integration of the KP hierarchy in a non formal way. One of the main advances of this work is to get a (non formal) principal bundle where the concept of holonomy makes sense rigorously. The geometric objects under consideration are diffeological or Fr\"olicher groups, which are regular in the sense that the exponential map exists and is smooth. 

Diffeological spaces, first described in the 80's by Souriau and his coworkers \cite{Don, Igdiff, Les, Sou} are generalizations of manifolds that enables differential geometry without charts. Independently, Fr\"olicher spaces give a more rigid framework, that also generalize the notion of manifolds \cite{CN, FK, KM}. The comparison of the two frameworks has been made independently in \cite{Ma2006-3} and in \cite{Wa}, see e.g. \cite{Ma2013}. 

We summarize and enlarge the theoretical results of \cite{Ma2013}, secondly giving three examples:

- one that enlarges straightway the example of \cite{Ma2013}: in $q-$deformation, indexes are monomials $q^n$, and in the example that we develop the base algebra is $A=Cl(M,E),$ the algebra of classical pseudo-differential operators. It is not an enlargeable Lie algebra, but the formal series $$G = 1+\sum_{n >0} q^n A$$ is a regular Fr\"olicher Lie group (section \ref{opd}).

- But one can also  consider indexes that are obtained replacing $\N$ by a groupo\"id: this is the case of path-like or cobordism-like formal series (section \ref{pc} and after). 
In this example, the strategy is based on defining a $\N-$grading on the groupo\"id, such that only the neutral element is of order $0.$ On such a groupo\"id, no element is invertible except the neutral element. This is the case for the cobordism composition, when we work on well-chosen families of manifolds $M,$ viewed as morphisms of cobordism, 
that are embedded in a fixed target space $N.$  The grading is induced by the volume, and we have to choose an algebra $A.$ In order to show that this setting is not void, 
we start from works \cite{AHK,AHKH1,AHKH2,AHKH3,AHKH4} that introduced so-called stochastic cosurfaces, that have applications in fields of mathematical physics such as 2D-yang-Mills theory, lattice models and Higgs fields among others. We enlarge the settings of the previous references, and adapt them to build families of measures indexed by cobordism, such that, if $\gamma $ and $\gamma'$ are two morphisms of cobordism that can be composed into $\gamma \gamma',$ then for the corresponding measures, we get $\mu(\gamma\gamma') = \mu(\gamma)\mu(\gamma')$ (convolution product). By the way, the mapping $\gamma \mapsto \mu(\gamma)$ can be understood as a formal serie over a family $\Gamma$ of morphisms for cobordisms. 

\vskip 12pt
These two examples have been developped together here because the technical tools necessary for exposition are very similar, and because we can already announce work in progress on applications of what we have called the dimension extension of the cosurface holonomy, linked with higher dimensional field theories, and independently, reformulation and extension of results of renormalization theory using the first example.
We finish this introduction by quoting two works \cite{BS2014,BDS2015} on the Butcher group and some of its generalizations, that are very similar to what we develop in this article. Initially, the author was not aware of these papers and intended to make also remarks on the Butcher group. But these two works are more complete, and the section on the Butcher group has been deleted from the text.

\section*{Acknowledgements} I would like to thank Professor Ambar Sengupta for stimulating discussions on the topics of cobordism, that influenced the corresponding section of this paper.
These discussions mostly occured during two stays at the Hausdorff Center f\"ur Mathematik at Bonn, Germany, invited by Sergio Albeverio and Matthias Lesch who are warmly acknowledged.    
	
\section{Regular Fr\"olicher Lie groups of series of unbounded operators} 

\subsection{Diffeological spaces and
Fr\"olicher spaces}
\label{1.1}

\begin{Definition} Let $X$ be a set.

\noindent $\bullet$ A \textbf{plot} of dimension $p$ (or $p$-plot)
on $X$ is a map from an open subset $O$ of $\R^{p}$ to $X$.

\noindent $\bullet$ A \textbf{diffeology} on $X$ is a set $\p$
of plots on $X$ such that, for all $p\in\N$,

- any constant map $\R^{p}\rightarrow X$ is in $\p$;

- Let $I$ be an arbitrary set; let $\{f_{i}:O_{i}\rightarrow X\}_{i\in I}$
be a family of maps that extend to a map $f:\bigcup_{i\in I}O_{i}\rightarrow X$.
If $\{f_{i}:O_{i}\rightarrow X\}_{i\in I}\subset\p$, then $f\in\p$.

- (chain rule) Let $f\in\p$, defined on $O\subset\R^{p}$. Let $q\in\N$,
$O'$ an open subset of $\R^{q}$ and $g$ a smooth map (in the usual
sense) from $O'$ to $O$. Then, $f\circ g\in\p$.

\vskip 6pt $\bullet$ If $\p$ is a diffeology $X$, $(X,\p)$ is
called \textbf{diffeological space}.

\noindent Let $(X,\p)$ et $(X',\p')$ be two diffeological spaces,
a map $f:X\rightarrow X'$ is \textbf{differentiable} (=smooth) if
and only if $f\circ\p\subset\p'$. \end{Definition}

We now introduce Fr\"olicher spaces.

\begin{Definition} $\bullet$ A \textbf{Fr\"olicher} space is a triple
$(X,\F,\C)$ such that

- $\C$ is a set of paths $\R\rightarrow X$,

- A function $f:X\rightarrow\R$ is in $\F$ if and only if for any
$c\in\C$, $f\circ c\in C^{\infty}(\R,\R)$;

- A path $c:\R\rightarrow X$ is in $\C$ (i.e. is a \textbf{contour})
if and only if for any $f\in\F$, $f\circ c\in C^{\infty}(\R,\R)$.

\vskip 5pt $\bullet$ Let $(X,\F,\C)$ et $(X',\F',\C')$ be two
Fr\"olicher spaces, a map $f:X\rightarrow X'$ is \textbf{differentiable}
(=smooth) if and only if $\F'\circ f\circ\C\in C^{\infty}(\R,\R)$.
\end{Definition}

Any family of maps $\F_{g}$ from $X$ to $\R$ generate a Fr\"olicher
structure $(X,\F,\C)$, setting \cite{KM}:

- $\C=\{c:\R\rightarrow X\hbox{ such that }\F_{g}\circ c\subset C^{\infty}(\R,\R)\}$

- $\F=\{f:X\rightarrow\R\hbox{ such that }f\circ\C\subset C^{\infty}(\R,\R)\}.$

This notion will be useful
in the sequel to describe in a simple way a Fr\"olicher structure.
A Fr\"olicher space, as a differential space, carries a natural topology,
which is the pull-back topology of $\R$ via $\F$. In the case of
a finite dimensional differentiable manifold, the underlying topology
of the Fr\"olicher structure is the same as the manifold topology. In
the infinite dimensional case, these two topologies differ very often.
One can remark, if $X$ is a Fr\"olicher space, we define a natural
diffeology on $X$ by \cite{Ma2006-3}, see e.g. \cite{Ma2013, Wa}: 
$$
\p(\F)=
\coprod_{p\in\N}\{\, f\hbox{ p-
paramatrization on } X; \, \F \circ f \in C^\infty(O,\R) \quad \hbox{(in
the usual sense)}\}.$$
With this construction, we also get a natural diffeology when
$X$ is a Fr\"olicher space. In this case, one can easily show the following:
\begin{Proposition}\cite{Ma2006-3}, see e.g. \cite{Ma2013, Wa}
Let$(X,\F,\C)$
and $(X',\F',\C')$ be two Fr\"olicher spaces. A map $f:X\rightarrow X'$
is smooth in the sense of Fr\"olicher if and only if it is smooth for
the underlying diffeologies. \end{Proposition}

Thus, we can also state:

\begin{tabular}{ccccc}
smooth manifold  & $\Rightarrow$  & Fr\"olicher space  & $\Rightarrow$  & Diffeological space\tabularnewline
\end{tabular}
\subsection{Fr\"olicher completion of a diffeological space} \label{complet}

We now finish the comparison of the notions of diffeological and Fr\"olicher 
space following mostly \cite{Wa}:

\begin{Theorem}
Let $(X,\p)$ be a diffeological space. There exists a unique Fr\"olicher structure
$(X, \F_\p, \C_\p)$ on $X$ such that for any Fr\"olicher structure $(X,\F,\C)$ on $X,$ these two equivalent conditions are fulfilled:

(i)  the canonical inclusion is smooth in the sense of Fr\"olicher $(X, \F_\p, \C_\p) \rightarrow (X, \F, \C)$

(ii) the canonical inclusion is smooth in the sense of diffeologies $(X,\p) \rightarrow (X, \p(\F)).$ 

\noindent Moreover, $\F_\p$ is generated by the family 
$$\F_0=\lbrace f : X \rightarrow \R \hbox{ smooth for the 
usual diffeology of } \R \rbrace.$$
\end{Theorem}

\noindent
\textbf{Proof.}
Let $(X,\F,\C)$ be a Fr\"olicher structure satisfying \textit{(ii)}. 
Let $p\in P$ of domain $O$. $\F \circ p \in C^\infty(O,\R)$ in the usual sense.  
Hence, if $(X,\F_\p, \C_\p)$is the Fr\"olicher structure on $X$ generated by the 
set of smooth maps $(X,\p)\rightarrow \R,$ we have two smooth inclusions 
$$ (X,\p) \rightarrow (X,\p(\F_\p)) \hbox{ in the sense of diffeologies }$$
and then
$$ (X, \F_\p, \C_\p) \rightarrow (X,\F,\C) \hbox{ in the sense of Fr\"olicher. }$$
 \qed

\begin{Definition} \cite{Wa}
A \textbf{reflexive} diffeological space is a diffeological space $(X,\p)$ such that $\p = \p(\F_\p).$
\end{Definition}

\begin{Theorem} \cite{Wa}
The category of Fr\"olicher spaces is exactly the category of reflexive diffeological spaces.
\end{Theorem}

This last theorem allows us to make no difference between Fr\"olicher spaces and reflexive diffeological spaces. 
We shall call them Fr\"olicher spaces, even when working with their underlying diffeologies.

\subsection{Push-forward, quotient and trace}

We give here only the results that will be used in the sequel. For an overview on diffeologies, see \cite{Sou} or more recently \cite{Igdiff}

\begin{Proposition} \cite{Ma2006-3} Let $(X,\p)$ be a diffeological space,
and let $X'$ be a set. Let $f:X\rightarrow X'$ be a surjective map.
Then, the set \[
f(\p)=\{u\hbox{ such that }u\hbox{ restricts to some maps of the type }f\circ p;p\in\p\}\]
 is a diffeology on $X'$, called the \textbf{push-forward diffeology}
on $X'$ by $f$. \end{Proposition}

We have now the tools needed to describe the diffeology on a quotient:

\begin{Proposition} \label{quotient} let $(X,\p)$ b a diffeological
space and $\rel$ an equivalence relation on $X$. Then, there is
a natural diffeology on $X/\rel$, noted by $\p/\rel$, defined as
the push-forward diffeology on $X/\rel$ by the quotient projection
$X\rightarrow X/\rel$. \end{Proposition}

Given a subset $X_{0}\subset X$, where $X$ is a Fr\"olicher space
or a diffeological space, we can define on trace structure on $X_{0}$,
induced by $X$.

$\bullet$ If $X$ is equipped with a diffeology $\p$, we can define
a diffeology $\p_{0}$ on $X_{0}$ setting \[
\p_{0}=\lbrace p\in\p\hbox{such that the image of }p\hbox{ is a subset of }X_{0}\rbrace.\]

$\bullet$ If $(X,\F,\C)$ is a Fr\"olicher space, we take as a generating
set of maps $\F_{g}$ on $X_{0}$ the restrictions of the maps $f\in\F$.
In that case, the contours (resp. the induced diffeology) on $X_{0}$
are the contours (resp. the plots) on $X$ which image is a subset
of $X_{0}$.

\subsection{Cartesian products and projective limits}

The category of Sikorski differential spaces is not cartesianly closed,
see e.g. \cite{CN}. This is why we prefer to avoid the questions
related to cartesian products on differential spaces in this text,
and focuse on Fr\"olicher and diffeological spaces, since the cartesian
product is a tool essential for the definition of configuration spaces.

In the case of diffeological spaces, we have the following \cite{Sou}:

\begin{Proposition} \label{prod1} Let $(X,\p)$ and $(X',\p')$
be two diffeological spaces. We call \textbf{product diffeology} on
$X\times X'$ the diffeology $\p\times\p'$ made of plots $g:O\rightarrow X\times X'$
that decompose as $g=f\times f'$, where $f:O\rightarrow X\in\p$
and $f':O\rightarrow X'\in\p'$. \end{Proposition}

Then, in the case of a Fr\"olicher space, we derive very easily, compare
with e.g. \cite{KM}: \begin{Proposition} \label{prod2} Let $(X,\F,\C)$
and $(X',\F',\C')$ be two Fr\"olicher spaces, with natural diffeologies
$\p$ and $\p'$ . There is a natural structure of Fr\"olicher space
on $X\times X'$ which contours $\C\times\C'$ are the 1-plots of
$\p\times\p'$. \end{Proposition}

We can even state the same results in the case of infinite products,
in a very trivial way by taking the cartesian products 
of the plots or of the contours.  Let us now give the description of what happens
for projective limits of Fr\"olicher and diffeological spaces.





\subsection{Regular Fr\"olicher groups of series}

Let $(G, \F, \C)$ be a Fr\"'olicher space which is a group  such that
the group law and the inversion map are smooth. These laws are also smooth
for the underlying diffeology. Then, following \cite{Les}, this is
possible as in the case of manifolds to define a tangent space and
a Lie algebra $\mathfrak{g}$ of $G$ using germs of smooth maps.
Let us precise the algebraic, diffeological and Fr\"olicher structures of $\mathfrak{g}.$

\begin{Proposition}
Let $\mathfrak{g} = \{ \partial_t c(0) ; c \in \C \hbox{ and } c(0)=e_G \}$
be the space of germs of paths at $e_G.$
\begin{itemize}
	\item Let $(X,Y) \in \mathfrak{g}^2,$ $X+Y = \partial_t(c.d)(0)$  where $c,d \in \C ^2,$ $c(0) = d(0) =e_G ,$ 
	$X = \partial_t c(0)$ and $Y = \partial_t d(0).$ 
	\item Let $(X,g) \in \mathfrak{g}\times G,$ $Ad_g(X) = \partial_t(g c g^{-1})(0)$  where $c \in \C ,$ $c(0) =e_G ,$ 
	and $X = \partial_t c(0).$  
	\item Let $(X,Y) \in \mathfrak{g}\times G,$ $[X,Y] = \partial_t( Ad_{c(t)}Y)$   where $c \in \C ,$ $c(0) =e_G ,$ 
	$X = \partial_t c(0).$
\end{itemize}
All these operations are smooth and thus well-defined. 
\end{Proposition} 

The basic properties remain globally the same as in the case of Lie groups, and the prrofs are similar replacing charts by plots of the underlying diffeologies. (see e.g. \cite{Les} for further details) 

\begin{Definition} A Fr\"olicher group $G$ with Lie algebra $\mathfrak{g}$
is called \textbf{regular} if and only if there is a smooth map \[
Exp:C^{\infty}([0;1],\mathfrak{g})\rightarrow C^{\infty}([0,1],G)\]
 such that $g(t)=Exp(v(t))$ if and only if $g$ is the unique solution
of the differential equation \[
\left\{ \begin{array}{l}
g(0)=e\\
\frac{dg(t)}{dt}g(t)^{-1}=v(t)\end{array}\right.\]
 We define \begin{eqnarray*}
exp:\mathfrak{g} & \rightarrow & G\\
v & \mapsto & exp(v)=g(1)\end{eqnarray*}
 where $g$ is the image by $Exp$ of the constant path $v.$ \end{Definition}

The classical setting for infinite dimensional differential geometry requires the model topological vector space
to be complete or Mac-Key complete. One of the reasons for this choice 
is to ensure the existence of the integral of a path over a 
compact interval. This means that the choice of an adaquate topology is 
necessary. For vector spaces, the basis for such a study can be found 
in \cite{KM}, when the properties of the so-called ``convenient vector spaces''
are given. We have to remark that a vector space for which addition and scalar multiplication
are compatible with a given Fr\"olicher structure needs only 
a topological structure to become a convenient vector space.
In order to circumvent these topological considerations, and adapting the     
terminology of regular Lie groups to vector spaces (which are viewed as abelian Lie groups), we set:

\begin{Definition}
Let $(V,\F, \C)$ be a Fr\"olicher vector space, i.e. a vector space $V$ equipped with a Fr\"'olicher structure compatible
with the vector space addition and the scalar multiplication. $(V,\F, \C)$ is \textbf{regular} if the Riemann integral is a smooth map 
$$ \int_0^{(.)} : C^\infty([0;1];V) \rightarrow C^\infty([0;1],V)$$ and satisfies the property:

 $\int_0^{(.)}v = u$ if and only if $u$ is the unique solution of 
the differential equation
\[
\left\{ \begin{array}{l}
u(0)=0\\
u'(t)=v(t)\end{array}\right. .\]

\end{Definition}

This definition is of course fulfilled if $V$ is a complete locally convex topological vector space, equipped with its natural Fr\"olicher structure.

\begin{Definition}
Let $G$ be a Fr\"olicher Lie group with Lie algebra $\mathfrak{g}.$ Then, $G$ is \textbf{regular with regular Lie algebra}
if both $G$ and $\mathfrak{g}$ are regular.
\end{Definition}

The first known example is the following \cite{Ma2013}:
\begin{Proposition} \label{omo} Let $(G_{n})_{n\in\N}$ be a sequence of Banach
Lie groups, increasing for $\supset,$ and such that the inclusions
are Lie group morphisms. Let $G=\bigcap_{n\in\N}G_{n}.$ Then, $G$
is a Fr\"olicher regular Lie group with regular Lie algebra $\mathfrak{g}=\bigcap_{n\in\N}\mathfrak{g}_{n}.$
\end{Proposition} 
Let us notice that there exists non regular Fr\"olicher Lie groups, see \cite{Ma2013-2}, where as there is no example of Fr\'echet Lie group that has been proved to be non regular \cite{KM}.
We now turn to  key results from \cite{Ma2013}: 

\begin{Theorem} \label{regulardeformation}
Let $(A_n)_{n \in \N^*} $ be a sequence of complete locally convex (Fr\"olicher)
vector spaces which are regular, 
equipped with a graded smooth multiplication operation
on $ \bigoplus_{n \in \N^*} A_n ,$ i.e. a multiplication such that 
$A_n .A_m \subset A_{n+m},$ smooth with respect to the corresponding Fr\"olicher structures.
Let us assume that:

Then, the set 
$$1 + \A = \left\{ 1 + \sum_{n \in \N^*} a_n | \forall n \in \N^* , a_n \in A_n \right\} $$
is a Fr\"olicher Lie group, with regular  Fr\"olicher Lie algebra
$$\A= \left\{ \sum_{n \in \N^*} a_n | \forall n \in \N^* , a_n \in A_n \right\}.$$
Moreover, the exponential map defines a bijection $\A \rightarrow 1+\A.$  
\end{Theorem}
\begin{Theorem}\label{exactsequence}
Let 
$$ 1 \rightarrow K \underrightarrow{i} G \underrightarrow{p}  H \rightarrow 1 $$
be an exact sequence of Fr\"olicher Lie groups, such that there is a smooth section $s : H \rightarrow G,$ and such that 
the trace diffeology from $G$ on $i(K)$ coindides with the push-forward diffeology from $K$ to $i(K).$
We consider also the corresponding sequence of Lie algebras
$$ 0 \rightarrow \mathfrak{k} \underrightarrow{i'} \mathfrak{g} \underrightarrow{p}  \mathfrak{h} \rightarrow 0 . $$
Then, 
\begin{itemize}
\item The Lie algebras $\mathfrak{k}$ and $\mathfrak{h}$ are regular if and only if the
Lie algebra $\mathfrak{g}$ is regular;
\item The Fr\"olicher Lie groups $K$ and $H$ are regular if and only if the Fr\"olicher Lie group $G$ is regular.
\end{itemize}

\end{Theorem}

\subsection{Groups of series that are regular Fr\'echet Lie groups}
We now asume that the algebras $\A_n$ of Theorem \ref{regulardeformation} are Fr\'echet vector spaces, 
and that the bilinear multiplication $$\A_n \times \A_m \rightarrow \A_{n+m}$$ is smooth.
Then we get the following theorem: 

\begin{Theorem} \label{regulardeformation2}
The group $1 + \A$ is a regular Fr\'echet Lie group with Lie algebra $\A.$ 
\end{Theorem}

\noindent
\textbf{Proof.} The exponential $\A \rightarrow 1+\A$ is already shown to be bijective, and $\A$, seen as a vector space, is endowed with the semi-norms on each $\A_n.$
Let us now show that the exponential of paths $C^\infty([0;1],\A) \rightarrow C^\infty([0;1];1+\A)$ is smooth, which will complete the proof.
Let $v \in C^\infty([0,1];\A).$  Let $s \in [0;1]$ and let $j = \lfloor ns \rfloor.$ We define
$$u_n(s) = \left(1+\left(s - \frac{j}{n}\right)v\left(\frac{j}{n}\right)\right) 
\prod_{i = 1}^{j}  \left(1+\frac{1}{n} v\left(\frac{j-i}{n}\right)\right).$$

We have that 
$$\lim_{n \rightarrow +\infty}\partial_s u_n(s) . u_n^{-1}(s) = \lim_{n \rightarrow +\infty} v\left(\frac{j}{n}\right)\left(1+\left(s - \frac{j}{n}\right)v\left(\frac{j}{n}\right)\right)^{-1} = v(s). $$
Moreover, the $A_m$ component of the product converges to a sum of integrals of the type  
$$\int_{1\geq s_1 \geq ... \geq s_k\geq 0} \left[ \prod_{i = 1}^k v(s_i)\right]_m (ds)^k
$$
for $k \leq m,$ which shows the convergence to a path $u\in C^\infty([0;1];1+\A)$ satisfying 
$$ \partial_s u(s).u^{-1}(s) = v(s)$$ which smoothly depends on the path $v \in C^\infty([0;1],\A)$ in the Fr\'echet sense. \qed

\subsection{Examples of $q-$deformed pseudo-differential operators} \label{opd}
In our work of Lax-type equations, we use the following group from \cite{Ma2013}:Let $E$ be a smooth vector bundle over a
compact manifold without boundary M. We denote by  $ Cl(M, E) $ (resp.  $ Cl^k (M, E)
$) the space of
 classical pseudo-differential operators (resp.
classical pseudo-differential operators of order k) acting on smooth
sections of $E$. We denote by $Cl^*(M,\mathbb{C}^n)$,
 $Cl^{0,*}(M,\mathbb{C}^n)$ the groups of
the units of the algebras $Cl(M,\mathbb{C}^n)$ and
$Cl^{0}(M,\mathbb{C}^n)$.
Notice that $Cl^{0,*}(M,\mathbb{C}^n)$ is a CBH Lie group, and belong 
to a wider class of such groups that is studied in \cite{Glo}.

.
\begin{Definition}
Let $q$ be a formal parameter. 
We define the algebra of formal series 
$$Cl_q(M,E) = \left\{ \sum_{t \in \N^*} q^k a_k | \forall k \in \N^*, a_k \in Cl(M,E) \right\}.$$
\end{Definition}
This is obviously an algebra, graded by the order (the valuation) into the variable  $q.$ Thus, setting
$$ \A_n = \left\{ q^n a_n | a_n \in Cl(M,E)\right\} ,$$
we can set $\A = Cl_q(M,E)$ and state the following consequence of Theorem \ref{regulardeformation2}:

\begin{Cor}
The group $1 + Cl_q(M,E)$ is a regular Fr\'echet Lie group with Lie algebra $Cl_q(M,E).$
\end{Cor}
Let $Cl^{0,*}(M,E)$ be the Lie group of invertible pseudo-differential operators of order 0. This group is known to be a regular Lie group since Omori, but the most efficient proof is actually in \cite{Glo}, to our knowledge.
We remark a short exact sequence of Fr\"olicher Lie groups:
$$ 0 \rightarrow 1 + Cl_q(M,E) \rightarrow Cl^{0,*}(M,E) + Cl_t(M,E) \rightarrow Cl^{0,*}(M,E) \rightarrow 0,$$  
which satisfies the conditions of Theorem \ref{extension}, in its version for regular Lie groups that can be found in \cite{KM}. Thus, we have the following:
\begin{Theorem}
$Cl^{0,*}(M,E) + Cl_q(M,E)$ is a regular Lie group with Lie algebra $Cl^{0}(M,E) + Cl_q(M,E).$
\end{Theorem}

\begin{rem} One could also develop a similar example, 
which could stand as a generalized version, 
with log-polyhomogeneous pseudo-differential operators or with other algebras of non classical operators, as desired. 
These examples are not developed here in order to avoid some too long lists of examples 
constructed in the same spirit.
\end{rem}

\section{Algebras and groups of series over groupoids}
We mimick and extend the procedure used in \cite{Ma2013}  
Let $(I,*)$ be a groupoid with neutral element $e.$
Let $A_i$ be a family of regular Fr\"olicher vector spaces indexed by $I.$ The family $\{\A_i; i \in I \}$ is equipped with a multiplication, associative and distributive with respect to addition in the vector spaces $\A_i,$ such that 
$$A_i . A_j \left\{ \begin{array}{ll} \subset A_{i*j} & \hbox{ if } i*j \hbox{ exists} \\
= 0 & \hbox{ otherwise. } \end{array} \right.$$and smooth for the Fr\"olicher structures.
Let $\mathcal{A}$ be the vector space of formal series of the type $$ a = \sum_{i \in I} a_i \quad ; \quad a_i \in \mathcal{A}_i $$
such that, for each $k \in I $, there is a finite number of indexes $(i,j) \in I^2$ such that $i*j=k .$
Notice that, with such a definition, $\A_e$ is an algebra.
From now, we assume $\A$ unital, and we note its unit element $1.$ We are not sure that, with this kind of definition, the exponential exists. In order to make the previous theorems valid, 
we have to define an adequate $\N-$grading.

\begin{Definition}
Let $I$ as above, such that, there is a $\N-$grading, that is,  a morphism of groupo\"ids $ord : I \rightarrow \N,$ such that $ord^{-1}(0) = \{e\}.$  
Let $A_i$ be family of regular Fr\"olicher vector spaces indexed by $I.$
Let $$\mathcal{A}_0 \subset \A = \left\{ \sum_{i \in I-\{e\}} a_i | a_i \in \mathcal{A}_i \right\}$$ be an algebra. The Fr\"olicher vector space $\mathcal{A}_0$ is called Fr\"olicher $I-$graded regular algebra if and only if it is equipped with a multiplication, associative and distributive with respect addition, smooth for the induced Fr\"olicher structure.
\end{Definition}

\begin{Proposition} \label{AI}
Let $\mathcal{A}$ be a Fr\"olicher $I-$graded regular algebra. It is the Lie algebra of the Fr\"olicher regular Lie group $1+\mathcal{A}.$
\end{Proposition}

\noindent
\textbf{Proof.}
This is a straight application of theorem \ref{regulardeformation}
\qed 

\begin{Theorem} \label{extension}
Let $\mathcal{A} = \bigoplus_{i \in I} \mathcal{A}_i$ be a Fr\"olicher $I-$graded regular algebra. Let $G$ be a regular Fr\"olicher Lie group, acting on $\mathcal{A}$ componentwise. Then,  $$G \oplus A$$ is a regular Fr\"olicher Lie group.
\end{Theorem}

\noindent
\textbf{Proof.}
Considering the exact sequence 
$$ 0 \rightarrow 1 + \A \rightarrow G \oplus \A \rightarrow G \rightarrow 0$$
there  is a (global) slice $G \rightarrow G \oplus \{0_\A \}$ so that Theorem \ref{exactsequence} applies . \qed
 
\begin{rem}\label{tens}

Notice that this small section is written in a heuristic way, in order to describe another kind of example that will be described rigorously in future works. Let $I$ be a family of manifolds, stable and finitely generated under cartesian product.  
Cartesian product is the composition law. We remark that it is graded with respect to the dimension of the manifold. A standard singleton can be added to $I$ as a neutral element of dimension $0$.
Let $J$ be a family of finite rank vector bundles over the family of manifolds $I,$ stable under tensor product. The scalar field can be understood as a neutral element of dimension $0.$
   By the way, Proposition \ref{AI} applies to the following algebras :
	
	- algebras of smooth sections of the finite rank vector bundles of $J$
	
	- algebras of operators acting on these sections. 
	
	When $I = \left\{(S^1)^n |n \in \N^* \right\},$ $J = \left\{(S^1)^n \times \mathbb{C} |n \in \N^* \right\},$
	and when the algebras under consideration are $Cl(S^1,\mathbb{C})^{\otimes n}$, we recognize a framework in the vincinity of the example given at the end of \cite{Ma2008}. 

\end{rem}

\section{Path-like and Cobordism-like deformations} \label{pc}
Let us now consider the set

$$ Gr = \coprod_{m \in \N^*} Gr_m$$
where $Gr_m$ is the set of m-dimensional connected oriented manifolds $M$, possibly with boundary, where the boundaries $\partial M$ are separated into two disconnected parts: the initial part $\alpha(M)$ and the final part $\beta(M).$
Then, we have a composition law $*$, called cobordism composition in the rest of the text, defined by the following relation:

\begin{Definition}
Let $m \in \N^*.$ Let $M, M' \in Gr_m.$ Then $M'' = M*M' \in Gr_m$ exists if 
\begin{enumerate}
\item $\alpha(M) = \beta(M') \neq \emptyset,$ up to diffeomorphism
\item $\alpha(M'')= \alpha(M')$
\item $\beta(M'') = \beta(M)$
\item $M''$ cuts into two pieces $M'' = M \cup M'$ with $M \cap M' = \alpha(M) = \beta(M').$
\end{enumerate}
\end{Definition} 

This composition, that we call \textbf{cobordism composition}, extends naturally to embedded manifolds:
\begin{Definition}
Let $N$ be a smooth manifold.
$$ Gr(N) = \coprod_{m \in \N^*} \coprod_{M \in Gr_m} Emb(M,N).$$
where the notation $Emb(M,N)$ denotes the smooth manifold of smooth embeddings of $M$ into $N.$
\end{Definition}
Notice that if $dim(N) < \infty,$ we have $m \leq dim(N).$
We recall that that $Gr(N)$ is naturally a smooth manifold, since $Emb(M,N)$ is a smooth manifold \cite{KM}, and that $*$ is obviously smooth because it is smooth in the sense of the underlying diffeologies. 
When we only consider manifolds without boundary (in this case,  cobordism composition is not defined), these spaces are called non linear grassmanians in the litterature, which explains the notations.
Let us now turn to $q-$deformed groups and algebras:
\begin{Definition}
\begin{itemize}
\item Let $ I =  \left( Gr \times \N^*\right) \coprod (\emptyset, 0),$ graded by the second component.  Assuming $\emptyset$ as a neutral element for $*$, we extend the cobordism composition into a composition, also noted $*$, defined as:
$$ (M,p) * (M',p') = (M*M', p+ p')$$ when $M*M'$ is defined. WE call \textbf{length} of $(M,p)$ the number $len(M,p)=p.$
\item
Let $ I(N) =  \left( Gr(N) \times \N^*\right) \coprod (\emptyset, 0),$ graded by the second component.  Assuming $\emptyset$ as a neutral element for $*$, we extend the cobordism composition into a composition, also noted $*$, defined as:
$$ (M,p) * (M',p') = (M*M', p+ p')$$ when $M*M'$ is defined.
\item Let $m \in \N^* .$ We note by $I_m$ and $I_m(N)$ the set of indexes based on $Gr_m$ and on $Gr_m(N)$ respectively
\end{itemize}
\end{Definition}
  Let $A$ be a regular Fr\"olicher algebra. Let $m \in \N^*.$ Let 
$$ \A_{I_m} = \left\{ \sum_{(M,n) \in I_m} q^na_{M,n} | a_{M;n} \in A \right\} $$
and let 
$$ \A_{I_m}(N) = \left\{ \sum_{(\phi,n) \in I_m(N)} q^na_{\phi,n} | a_{M;n} \in A \right\}. $$
\begin{Theorem} \label{cob}
Let $\Gamma \subset \coprod_{m \in \N^*} I_m$, resp. $\Gamma(N) \subset \coprod_{m \in \N^*} I_m(N)$, be a family of indexes, stable under $*,$ such that $\forall m \in \N^*,$
\begin{enumerate}
\item  $\forall m \in \N^*,$ $\Gamma \cap I_m$ is finite or, more generally;
\item $\forall \gamma \in \Gamma,$ the set of pairs $(\gamma',\gamma'') \in \Gamma^2$ such that $\gamma = \gamma' *`\gamma''$ is finite.
\item 
$$ \A_\Gamma = \left\{ \sum_{(M,n) \in \Gamma} q^na_{M,n} | a_{M;n} \in A \right\} ,$$ 
resp. $$ \A_\Gamma(N) = \left\{ \sum_{(\phi,n) \in \Gamma(N)} q^na_{\phi,n} | a_{\phi;n} \in A \right\} ,$$
is a regular algebra (in the sense of regular vector spaces). 
\end{enumerate}
  Then
 $$ 1_A + \A_{\Gamma-\{(\emptyset,0)\}}$$
is a Lie group with Lie algebra $A_\Gamma.$
  Moreover, for each Fr\"olicher Lie group $G $ with Lie algebra $\mathfrak{g}$ such that $G \subset A^*$ smoothly, $$ G \oplus \A_{\Gamma - \{(\emptyset, 0)\}}$$
is a regular Fr\"olicher Lie group with Lie algebra $\mathfrak{g} \oplus \mathcal{A}_{\Gamma-\{(\emptyset,0)\}}.$
Moreover, the results are the same replacing $\Gamma$ by $\Gamma(N).$
\end{Theorem} 	

\noindent
\textbf{Proof.}
Following condition (1) or (2), in the (possibly infinite sum) $\sum_{(M,n) \in \Gamma} q^na_{M,n},$ each power $q^n$ has only a finite number of $A-$coefficients since there is only a finite number of possible indexes for each $q^n.$ So that, Proposition \ref{AI} and Theorem \ref{extension} apply. The same arguments are also valid when replacing $\Gamma$ by $\Gamma(N).$ \qed



\section{Markov Cosurfaces in codimension 1}\label{s:cosurface}

 \subsection{Settings}
Let $M$ be a $d-$dimensional connected oriented
Riemannian manifold. Let $\mch_\vee$ be the set of
embedded, oriented,
smooth, closed, connected hypersurfaces (codimension $1$ submanifolds) of $M$
with piecewise smooth border.
What we call hypersurface is mostly smooth hypersurfaces
on the mnifold $M,$ but since we need piecewise smooth oriented hypersurfaces,
we need to build them by induction, gluing together the smooth components.
What we get at the end is a space of oriented piecewise smooth hypersurfaces, 
with piecewise smooth border.
\begin{Definition}\label{hyperabel}

We set $\mch_\vee^{(1)}=\mch_\vee.$
For $d\geq 2,$ we define by induction:

$\bullet$ Let $(s_1,s_2)\in \mch_\vee\times\mch_\vee.$
If
\begin{enumerate}
\item $s_1 \cap s_2 \subset \partial s_1 \cap \partial s_2$ is a
$(d-2)$ piecewise smooth manifold and

\item the orientations induced on $s_1 \cap s_2$ by $s_1$ and $s_2$
are opposite,
\end{enumerate}
then we define $s_1 \vee s_2$ to be the oriented
piecewise smooth hypersurface of $M$ obtained by
gluing $s_1$ and $s_2$ along their common border.
The orientation of $s_1\vee s_2$
is the one induced by $s_1$ and $s_2.$
The set of all such hypersurfaces is denoted by $\mch_\vee^{(2)}.$

 $\bullet$ Let $(s_1,s_2)\in \mch_\vee^{(n-1)}\times \mch_\vee.$
If
\begin{enumerate}
\item $s_1 \cap s_2 \subset \partial s_1 \cap \partial s_2$ is a
$(d-2)$ piecewise smooth manifold and

\item the orientations induced on $s_1 \cap s_2$ by $s_1$ and $s_2$
are opposite,
\end{enumerate}
then we define in the same way $s_1 \vee s_2.$
The set of such hypersurfaces is denoted by $\mch_\vee^{(n)}.$

$\bullet$ We set $\Sigma_\vee = \bigcup_{n \in \N^*}\mch_\vee^{(n)}.$
\end{Definition}

In all the article, we shall assume also that the connected components of 
$\partial s$ are in $\Sigma_{\vee,n-1}$ if $s \in \Sigma_{\vee,n},$ for $n \in \N^*.$
\vskip 12pt
\noindent

\begin{rem} If $M$ is $2-$dimensional, it might seem that definition \ref{hyperabel}
generalizes  the composition of unparametrized
piecewise smooth paths, setting $\h$ to be
the set of smooth paths and $\vee$
the groupo\"id composition law of
oriented piecewise smooth paths.
In fact, we need to reformulate
the definition for $d=2$ in order
to fit with the usual composition
of paths.
Let us look at the following example.
Let $M=\R^2$, and let $s_1$ and $s_2$
the paths parametrized
by $s_1(t) = (\cos(\pi t),\sin(\pi t))$ and
$s_2(t) = (-\cos(\pi t), -\sin(\pi t))$ for $t \in [0,1].$
We have $\partial s_1 = \partial s_2 =\{(-1;0);(1;0)\},$
with ``opposite orientations''
(i.e. the endpoint of $s_1$ (resp. $s_2$) is the
initial point of $s_2$ (resp. $s_1$))
so that paths $s_1$ and $s_2$ can be composed.
But in order to have a loop,
one has to determine which point among $\{(-1;0);(1;0)\}$
will be the initial point. The choice comes
with the order in the composition of paths: $s_1 * s_2$
or $s_2 * s_1.$ Such a choice cannot be done with
definition \ref{hyperabel}
because the law $\vee$ is obviously
commutative in $\mch_\vee \times\mch_\vee.$
\end{rem}
Now, we define
$\mch_*$ the set of (unparametrized, but oriented) smooth hypersurfaces $s$
on the oriented manifold $M,$ equipped in addition 
with a prescribed orientation of smooth components of
$\partial s .$
\textbf{Initial parts} of $\partial s$, noted $\alpha(s),$ are those for which the
prescribed orientation is opposite to the one induced by $s,$
and the \textbf{final parts}, noted $\beta(s),$
are the ones for which they co\"incide. For $d=2,$ the orientation
of paths can prescribe naturally initial and final points. This is the (apparently natural) choice 
that has been made in \cite{AHKH1} but we remark here that this choice
is not necessary. The picture of the following definition will be merely the same 
as the one of definition \ref{hyperabel}, but each smooth component of the 
border of the hypersurface is assigned to be either initial or final. In order 
to keep the coherence with the loop composition, 

- we can glue together a final part with an initial part,

- and the final parts and the initial parts can be the same set-theorically, 
just as in the case of a loop starting and finishing at the same point.   

Here is the construction:
\begin{Definition}\label{hypernonabel}

We set $\mch_*^{(1)}=\mch_*.$
we define by induction:

$\bullet$ Let $(s_1,s_2)\in \mch_*\times \mch_*.$
Let $a = \alpha(s_1)\cap \beta(s_2).$
We define $s_1 * s_2$ as the oriented
piecewise smooth hypersurface of $M$ obtained
by gluing $s_1$ and $s_2$ on $a $ and denoted by $s_1 \cup_a s_2.$
The orientation of $s_1 * s_2$
is the one induced by $s_1$ and $s_2$ on $s_1 \cup_a s_2$ 
By the way, we have $ \alpha(s_1 * s_2) = \alpha(s_2) \cup (\alpha(s_1) - a),$
$\beta(s_1 * s_2)= \beta(s_1) \cup (\beta(s_2)-a),$ and 
$\partial(s_1 * s_2) = \alpha(s_1*s_2) \coprod \beta(s_1 * s_2).$ 
The set of such hypersurfaces is denoted by $\mch_*^{(2)}.$

 $\bullet$ Let $(s_1,s_2)\in \mch_*^{(n-1)}\times \mch_*.$
then we define in the same way $s_1 * s_2.$
The set of such hypersurfaces is denoted by $\mch_*^{(n)}.$

$\bullet$ We set $\Sigma_* = \bigcup_{n \in \N^*}\mch_*^{(n)}.$
\end{Definition}

Notice that there is a forgetful map $\Sigma_* \rightarrow \Sigma_\vee$ 
only for the hypersurfaces $s \in \Sigma^*$ that have no self-intersection. 
The following example, based on the M\"obius band, shows how the pictures can 
turn badly if we forget this restriction.   
\noindent
\textbf{Example.}
Let us fix $M = \R^3$ and let $$s_1 = \left\{\left(\cos(\pi t), \sin(\pi t), s-\frac{1}{2}\right) | (t,s) \in [0;1]^2\right\}$$
such that $\alpha(s_1) = \left\{(\cos(\pi t), \sin(\pi t), s) | (t,s) \in (\partial[0;1])^2, t=0\right\},$
$\beta(s_1) = \partial s_1 - \alpha(s_1)$ and let 
 
$$s_2 = \left\{(-\cos(\pi t), -\sin(\pi t)(1+cos(\pi t)), s-cos(\pi t)) | (t,s) \in [0;1]^2\right\}$$
such that 
$$\alpha(s_2) = \left\{\left(-\cos(\pi t), -\sin(\pi t)\left(1+\frac{sin(\pi t)}{2}\right), s-cos(\pi t)-\frac{1}{2}\right) | (t,s) \in (\partial[0;1])^2, t=0\right\},$$
$\beta(s_2) = \partial s_2 - \alpha(s_2).$ If one glues topologically 
$s_1$ and $s_2,$ we get the M\"obius band which is non orientable. So that 
$s_1 \vee s_2$ is not defined in this case. By our choices of initial ad final parts,
$s_1 * s_2$ and $s_2 * s_1$ exist both, because they can be represented by 
the ``cut'' M\"obius band, with
$$\alpha(s_1 * s_2) = \alpha(s_2) = \left\{(1; 0; s)| 
-\frac{1}{2}\leq s \leq \frac{1}{2} \right\},$$ 
$$\beta(s_1 * s_2) = \beta(s_1) = \left\{(1; 0; s)| 
-\frac{1}{2}\leq s \leq \frac{1}{2} \right\},$$
 and with  
$$\alpha(s_2 * s_1) = \alpha(s_1) = \left\{(-1; 0; s)| -\frac{1}{2}\leq s \leq \frac{1}{2} \right\},$$ 
$$\beta(s_2 * s_1) = \beta(s_2) = \left\{(-1; 0; s)| -\frac{1}{2}\leq s \leq \frac{1}{2} \right\}.$$
One can say that this is not natural since e.g. $\alpha(s_1*s_2)$ is not in the border of the underlying $C^0
$-manifold, this is one of the reasons why we discuss this example in details. 
Moreover, this fits with the natural composition of paths: ignoring the third coordinate, we get back the 
classical composition of paths, for which loops are topologically without border but 
have a starting point and an endpoint.

\vskip 12pt 

In what follow, $\Sigma_M$ represents either $\Sigma_\vee$ or $\Sigma_*$
with an adequate choice of Lie group $G$ (we choose $G$ to be abelian for $\Sigma_\vee $).
Anyway, we denote the group law of $G$ by the operation of multiplication, and we note $s_1s_2$ for $s_1 \vee s_2$ or $s_1 * s_2.$
\begin{Definition}

A  $G-$valued \textbf{cosurface} is a map $$c : \Sigma_M \rightarrow G$$
such that
\begin{enumerate}
\item $\forall (s_1,s_2)\in \Sigma_M \times \h, c(s_1 \vee s_2) =c(s_1)c(s_2)$ and
\item We denote by $\tilde s$ the same hypersurface as
$s\in \Sigma_M$ with opposite orientation. Then
$\forall s \in \Sigma_M, c(\tilde s) = c(s)^{-1}.$
\end{enumerate}
Let $\tau_s(c)=c(s).$
Let $\Gamma_{M,G}$ the set of $G-$valued cosurfaces of $M$
equipped with the the smallest $\sigma-$ algebra making measurable
the collection of maps
$$\{\tau_s : \Gamma_{M,G} \rightarrow G |s \in \Sigma_M \}.$$
\end{Definition}
\vskip 12pt
Let $(\Omega, \mathcal{B}, p)$ be any probability space.

\begin{Definition}
A \textbf{stochastic cosurface} is a map
$$ C :\Omega \times \Sigma_M \rightarrow G $$
such that:
\begin{enumerate}
\item $\forall \omega \in \Omega,$ $C(\omega,.) \in \Gamma_{M,G}.$
\item the map $\omega \in \Omega \mapsto C(\omega,.)$
is a $\Gamma_{M,G}-$valued measurable map.
\end{enumerate}
\end{Definition}
For a subset $\Lambda \subset M$ we consider the $\sigma-$algebra $\mathfrak{T}(\Lambda)$ generated by stochastic cosurfaces $C(s)$ where $s \subset \Lambda.$
In other words,
$$\mathfrak{T}(\Lambda) = \sigma\left\{ \left\{C \in \Gamma_{M,G} | C(s)\in B\right\} | s\subset \Lambda; B \hbox{ Borel subset of } G \right\}.$$

Now, we have to define finite sequences of hypersurfaces that we denote by \textbf{complex}.

\begin{Definition}
Let $n \in \mathbb{N}^*.$ An $n-$\textbf{complex} on $M$ is a $n-$uple $K=(s_1,...s_n) \in \left(\Sigma_M\right)^n$ such that $s_i \neq s_j$ for $i\neq j.$ We define
$$ C(K) = (C_1(s_1),...,C_n(s_n))$$ where each $C_i \in \Gamma(M;G).$ We note
by $\mathcal{K}$ the set of complexes of any length $n.$
\end{Definition}

Notice that a complex is an \textbf{ordered} sequence, and related with this order
 there is a natural notion of \textbf{subcomplex} of a complex $K.$ 
If $K=(s_1,...,s_n),$ a subcomplex $L$ is a subsequence of $K,$ 
that is $$\exists l< m\leq n, L = (s_l,...,s_m)=(s_i)_{l\leq i \leq m}.$$

Now, we need to recognize the complexes that define skeletons of a 
partition of the manifold $M.$

\begin{Definition}
A $n-$complex $K = (s_1,...,s_n)$ is \textbf{regular} if
  $\forall (i,j) \in {\mathbb{N}_n}^2,$ $$ i \neq j \Rightarrow s_i \cap s_j \subset \partial s_i \cap \partial s_j.$$
\end{Definition}
Notice that the definition does not consider initial and final parts of the borders.
In the sequel, since many ways to understand complexes can be useful 
(set-theory, topological spaces, sequences, oriented manifolds) 
we shall use the standard notations in these various fields and we shall precise 
in what sense we use them only if the notations carry any ambiguity.

\begin{Definition}
Let $K$ be a regular $n-$complex. $K$ is called \textbf{saturated}
if and only if $\bigcup_{i \in \mathbb{N}_n}s_i$ defines the borders
of a covering of $M$ by connected and simply connected closed subsets.
In other words, there is a familly $(A_k)_k$ of closed connected and simply connected subsets of $M$
such that
\begin{enumerate}
\item $\bigcup_{k}A_k = M$
\item for two any indexes $k$ and $k',$
if $k \neq k',$ $$A_k \cap A_{k'} \subset \partial A_k \cap \partial A_{k'} \subset\bigcup_{i \in \mathbb{N}_n}s_i.$$
\end{enumerate}
\end{Definition}

We say that a regular $n-$complex $K$ \textbf{splits} 
$M$ through the subcomplex $L = (s_i)_{l \leq i \leq m}$
if
\begin{enumerate}
\item $\bigcup_{s \in L} s $ splits $M$ into two connected components $M^+$ and $M^-$ and

\item $\bigcup_{i<l}s_i\subset Adh(M^-) ;$ we note ${K^-} = (s_i)_{i<l} ``= K \cap {M^-} ",$

\item $\bigcup_{i>m}s_i\subset Adh(M^+);$ we note ${K^+} = (s_i)_{i>m} ``= K \cap {M^+} ".$
  \end{enumerate}
(here, $Adh$ means topological closure)

\vskip 12pt
\noindent
\textbf{Examples around the $2$-cube} $ABCDEFGH \subset \R^3$ {\bf :}

$ABCDEFGH$ is the $2$-cube, and we assume that each coordinate of 
$A,B,C,D,E,F,G$ and $H$ is equal to $\pm 1.$

- Let us consider the (empty) $2-$cube $ABCDEFGH$ as a piecewise smooth 
hypersurface of $\R^3.$ By the orientation of $\R^3,$ and since 
the cube divides $\R^3$ into an inside part and an outside part, each face 
is oriented so that $ABCDEFGH\in \Sigma_\vee.$ The $2$-cube divides $\R^3$
into two parts, that we recognize as interior and exterior, 
which are connected an simply connected (but one is non contractible). 
So that, it splits $\R^3.$ We can also say that we have a regular $6-$complex
made of the faces of $ABCDEFGH$ (where we have to choose an order which is non canonical).     

-Let us now project $ABCDEFGH$ into $S^2$ radially. Then the segments
of $ABCDEFGH$ define a class of regular complexes on $S^2.$ The complex of the segments is not uniquely 
defined because of the order that we have to choose, and also because of the orientations of the segments 
that we have to choose. This will yield different possible 
splittings. For example,if we consider the the complex 
$K=(s_1,...s_{12})$ defined
by 
\begin{eqnarray*} K &=& \big((A;B),(B;E),(B;C),(C;D),(C;G),(A;E),(E;F),\\
&&(F;G),(G;H),(H;D),(D;A),(E;H)
\big)
\end{eqnarray*}
we have the subcomplex 
\begin{eqnarray*} L &=& \left((A;E),(E;F),(F;G),(G;H),(H;D),(D;A)\right)
\end{eqnarray*}
that splits $K.$
We have $(S^2)^+ = AEHD \cup EFGH$ and $(S^2)^-=ABCD \cup ABFE \cup BCGF \cup CDGH$
(which are both contractible) and finally 
$K^+ = (E;H)$ and $K^-=\left((A;B),(B;E),(B;C),(C;D),(C;G)\right).$
\begin{Definition}

Let $C $ be a stochastic cosurface. $C$ is said to be 
a \textbf{Markov cosurface}
if for each $n \in \N^*$ and for each regular $n-$complex 
$K$ which splits through a subcomplex $L=(s_i)_{l \leq i \leq m},$ for 
each couple of map $(f^+,f^-)$ for which the above expectations exists 
and $f^+$ (resp. $f^-$) is $\mathfrak{T}(M^+ \cup \bigcup_{l \leq i 
\leq m} s_i)-$measurable (resp.  $\mathfrak{T}(M^- \cup \bigcup_{l \leq i \leq m} 
s_i)-$measurable), we have:
$$\mathbb{E}(f^+f^-|\mathcal{T}(\bigcup_{l \leq i \leq m} s_i)) = \mathbb{E}(f^+|\mathfrak{T}(\bigcup_{l \leq i \leq m} s_i))\mathbb{E}(f^-|\mathfrak{T}(\bigcup_{l \leq i \leq m} s_i)).$$
\end{Definition}
\subsection{Markov cosurfaces and Markov semigroups}
Let $\lambda$ be the Haar measure on $G$ which is now assumed unimodular.
We introduce a projective system of probability measures on $\{ G^{K}; K \in \mathcal{K}\}.$ For this, we use a partial order on $\mathcal{K}.$

\begin{Proposition}
Let $(K,K')\in \mathcal{K}^2$ such that $K=K'$ in the set-theoric sense. We write $K \prec K'$ if $\forall s \in K,$
there exists a subcomplex $L'$ of $K'$ such that $s$ is the composition
of the elements of $L',$ ordered by indexes.
$\prec$ is an order on $\mathcal{K}.$
\end{Proposition}
\vskip 12pt
\noindent
\textbf{Proof.} Comparing this proposition with \cite{AHKH4}, we already have that 
$\prec$ is only a preorder. So that, we need only to check reflexivity. 
Let $s \in K.$ Taking $L'=\{s\},$ we get $K\prec K.$
Moreover, let $(K,K') \in \mathcal{K}^2,$ if $K \prec K'$ and $K'\prec K,$ 
$\forall s \in K, s\in K'$ and $\forall s' \in K', s \in K$ 
and hence K and K' have the same  hypersurfaces, indexed with respect to the same order.\qed

\vskip 12pt
We now recall the standard definition of filters for the order $\prec.$
\begin{Definition} A filter $P \subset \mathcal{K}$ is such that:
\begin{enumerate}
\item $\forall (K,K')\in P^2, \exists K'' \in P, \left( K \prec K'' \wedge K' \prec K''\right).$
\item $\left(\forall K \in \mathcal{K}, \exists K' \in P, K \prec K'\right) \Rightarrow K \in P.$
\end{enumerate}
\end{Definition}
Let $Q_t$ be a convolution semigroup of probability measures on $G$ with densities, i.e. $Q_t = q_t . \lambda$ satisfies
\begin{enumerate}
\item $Q_0 = \delta_e$ (Dirac measure at the unit element)
\item $\forall s,t \in (\R_+^*)^2,$ $\forall x \in G,$ $(q_t q_s)(x) = \int_Gq_s(xy^{-1})q_t(y)d\lambda(y)= q_{t+s}(x)$
\item $\lim_{t\rightarrow 0} Q_t =\delta_e$ weakly
\item $\forall (x,y) \in G^2, q_{(.)}(xy) = q_{(.)}(yx)$
\end{enumerate}
Now, we need to separate the exposition among the two approaches of cosurfaces, 
one on $\Sigma_\vee$ and the other on $\Sigma_*.$ 
In both cases, we fix $K \in P$ a regular saturated complex with associated 
domains $D=\{A_1,...A_m\}.$
Each $A\in D$ is oriented through the orientation of $M.$

$\bullet$ \underline{On $\Sigma_\vee.$}

We define
$$ \varphi_{A}(s) = \left\{ \begin{array}{cl} s &\hbox{ if } s\subset\partial A \hbox{ has the same orientation as } \partial A \\
\tilde{s} & \hbox{ if } s\subset\partial A \hbox{ has the opposite orientation from } \partial A \\
\emptyset & \hbox{ if } s \not\subset \partial A \end{array} \right. $$
We set $$\phi_A(C(K)) = \prod_{s \in K} C\circ\varphi_A(s).$$
(this product is a convolution product of measures)

$\bullet$ \underline{On $\Sigma_*.$}
Here, $K=(s_1,...,s_m)$ is ordered by indexation. Then, we have to work by induction to define $\varphi_A.$
\begin{itemize}
\item Let $s_j\in K$ be the first element in $K$ such that 
$s_j\subset \partial A.$ For $i<j,$
we set $\varphi_A(s_i)=\emptyset.$ Then 
we compare the orientation of $\partial A$ with the one of $s_j$ 
as in the case of $\Sigma_M^\vee$ to define $\varphi_A(s_j).$
\item Assume now that we have determined $\varphi_A$ till the index $j.$ Take $l$
to be the first index after $j$ such that $s_l \subset \partial A.$ As before, 
for $j<i<l,$ we set $\varphi_A(s_i) = \emptyset.$ First, compare the orientation 
of $s_l$ with the one of $\partial A$ and change $s_l$ into $\tilde{s_l}$
if necessary as before. Notice that final parts of $s_j$ and initial parts of $s_l$
are not considered here.
This enables anyway to define
 \begin{equation} \label{phiA} \phi_A(C(K)) =  C\circ\varphi_A(s_1)...C\circ\varphi_A(s_n) \end{equation}

\end{itemize}

Both in the case of $\Sigma_\vee$ and in the case of $\Sigma_*,$ we set

\begin{Definition} \label{muK}$$\mu_K^Q(C) = \prod_{i =1}^k q_{|A_i|}\left(\phi_{A_i}\circ C(K)\right).$$\end{Definition}

\begin{rem}
When $d=2$, changing the orientation of the path $s \in \mathcal{H}$ is the same 
as permuting its initial and its final points. Then, the procedure described for 
$\Sigma_*$ makes also final parts and initial parts co\"inciding.
\end{rem}

\subsection{Action of the symmetric group}
Looking at the definition \ref{muK} of $\mu_K,$ we easily see that
the value of $\mu_K$ is independent of the order
of the sequences $A=(A_1,...,A_k)$ since the group $G$ is unimodular. 
Unlikely, there is no invariance
under reordering $K$ in the non abelian case
(see the definition of $\phi_A$ in equation \ref{phiA}).
Here, $\phi_A$ depends on the order of the saturated complex $K = (s_1,...,s_n).$
So that, the action of the $n$-symmetric group $\mathfrak{G}_n$ on
indexes of $n$-saturated complexes
$$\left(\sigma, K=(s_1,...s_n)\right)\longmapsto  \sigma.K = \left( s_{\sigma(1)},...,s_{\sigma(n)} \right)$$
 generates an action $\sigma.\mu_K = \mu_{\sigma. K}.$
Setting $\mathcal{KS}$ to be the set of saturated complexes on $M,$ 
noting by $\mathfrak{G}_{\infty}$
the group of bijections on $\N^*,$ we get an action
$\mathfrak{G}_\infty \times \mathcal{KS} \rightarrow \mathcal{KS}$
in the following way:
completing $K = (s_1,...,s_n) \in \mathcal{KS}$ into
$\hat{K}=(s_1,...,s_n, \emptyset,...) \in \Sigma_*^{\N^*},$
a bijection $\sigma \in \mathfrak{G}_\infty$ on indexes gives a sequence
$\sigma.\hat{K}$ with only $n$ elements different from $\emptyset.$
We define $\sigma.K$ to be the $n$-saturated complex 
(indexed by $\N_n=\{1,2,...,n\}$) as the
collection $\{s_1,...,s_n\},$ ordered $\sigma. \hat{K}$ indexwise.

\subsection{Examples}
This selection of examples is based on earlier works \cite{AHK,AHKH1,AHKH2,AHKH3,AHKH4} 
where only $d=2$ examples on $\Sigma_*$ or examples on $\Sigma_\vee$ were considered.
An example on $\Sigma_*$ with $d=3$ will be given later because the tools 
needed have to  be much clarified.
\subsubsection{The $d=2$ holonomy cosurface.} \label{sss:holcosurface}
(\cite{AHKH1}, compare with the settings described in \cite{Levy2012}) 
Let $M$ be a $2-$manifold.
Let $\p(M)$be the space of piecewise oriented smooth paths, with canonical 
initial and final points (``canonical'' means induced by the path orientation).
By the way, open paths in $\p(M)$ can be identified with $\Sigma_\vee$ and
there is a map $\Sigma_{*} \rightarrow \p(M) $ which co\"incides with
the forgetful map $\Sigma_* \rightarrow \Sigma_\vee$ on open paths
and  with a map on loops that is changing the initial and final parts if necessary. 
Let $G$ be a Lie group and let $P = M\times G$ the trivial principal 
bundle over $M$ with structure group $G.$
Let $\theta$ be a connection on $P$ and we note by $Hol_\theta$ 
the holonomy mapping $\p(M) \rightarrow G$ for which the horizontal lift 
starts at $(\alpha(p),e_{G}).$
Let $s \in \Sigma_*$ that we identify with the 
corresponding path that we note also $s.$ We define the holonomy  cosurface $c$
by
  $$c(s) = Hol_\theta(\tilde s).$$ 
We need here to invert the orientation of the path because of 
the right action of the holonomy group on the principal bundle $P.$

\begin{rem} Since the definition of (non stochastic) cosurfaces
does not require any measure, one can take for $G$ any Lie group
on which the notion of horizontal lift of a path with respect to 
a connection is well defined. At this level, the construction is 
valid for any (finite dimensional) Lie group, but for any Banach Lie group, 
or regular Fr\'echet Lie group, regular $c_\infty$-Lie group 
\cite{KM}, as well as for 
regular fr\"olicher Lie group \cite{Ma2013}. 

\end{rem}  

\begin{rem}
	This approach is quite similar to the approach of gauge theories via quantum loop gravity approach, see e.g. \cite{Ma2018-3} for an up-to date paper. However, much open questions remain when one wish to work along the lines of this viewpoint.
\end{rem}
Then choosing $Q$ as the heat semi-group on $G,$
we get a stochastic cosurface picture of $2d$-Yang-Mills fields 
(see e.g. \cite{Levy2012}, \cite{Sengupta} and references therein 
for an extensive work when $M$ is a $2$-dimensional 
manifold and the topology is non trivial).

\subsubsection{Markov cosurfaces and lattice models} \cite{AHKH1,AHKH4}
Let $L_\epsilon = \Z^d.$ Let $U$ be an invariant function on a compact group 
$G$ and a ``coupling constant'' $\beta > 0.$ Let $\Lambda$ be a bounded subset 
of $L_\epsilon$ and let us define the (normalized) probability measure
$$ \mu_\epsilon^\Lambda = \frac{1}{Z_{\Lambda,\epsilon}} 
\exp\left(-\beta \sum_{\gamma \subset \Lambda} U(C(\partial\gamma))\right) 
\prod_{\gamma \subset \Lambda} dC(\partial \gamma)$$
where $\gamma$ is an elementary cell, $\partial \gamma$ the boundary 
$C(\partial \gamma)$ a variable associated to $\partial \gamma$ with values in $G.$
In the sens of projective limits of measures, 
the limit $\Lambda \rightarrow L_\epsilon$ exists and defines a 
Gibbs-like lattice cosurface. 
In the cases $G=U(1), SU(2),$ the continuum limit $\epsilon \rightarrow 0$ 
for $\mu_\epsilon$ has been shown to exist 
for appropriate $U$ and by a suitable choice of $\beta(\epsilon)$ such that 
$\lim_{\epsilon \rightarrow 0} \beta(\epsilon) = +\infty.$

\subsubsection{Markov cosurfaces and Higgs fields in 2-dimensional space time}(\cite{AHKH4}, see also \cite{AHK})
Let $\Lambda$ be a bounded non empty subset of $\mathbb{Z}^2.$
Cosurfaces $C$ are defined along the edges of $\mathbb{Z}^2.$
Let $G$ be a compact Lie group, equipped with a representation $\rho$ on an Euclide space $V$ with scalar product $<.;.>$ and norm $|.|.$
Let $(\lambda,\mu)\in \left(\mathbb{R}_+^*\right)^2$ and let
$\varphi$ be a $V-$ valued random field over $\mathbb{Z}^2.$
We define a probability measure$\mu_{\lambda,\Lambda}$ on $\mathbb{Z}^2$ with support in $\Lambda$
by
$$ \mu_{\lambda,\Lambda}(d\phi) = \frac{1}{Z_\Lambda}e^{-\frac{\lambda}{2}\sum_{x \in \Lambda} \left(B + \frac{\mu^2}{\lambda}\right)|\varphi(x)|^2} e^{-\frac{\lambda}{2}\sum_{x,y \in \Lambda}<\varphi(x),\rho \circ C(xy )\varphi(y)>}\prod_{x \in \Lambda} d\varphi(x).$$
where $C(xy)$ is the evaluation of the Markov cosurface $C$ on the path $xy.$
Replaing the lattice $\Z^2$ by $\epsilon \Z^2,$ with suitable choice of 
$\lambda(\epsilon, \mu(\epsilon),$ yields continuum limit Higgs models as 
$\epsilon \rightarrow 0.$ In the cases $G=U(1), SU(2),$ with $Q_t$ the heat 
semi-group, our Higgs fields coincide with the Higgs fields in the 
physics litterature.
For a description of the mathematical construction of the continuum limit, 
keeping $C$ fixed, see \cite{AHKH4}. 

\section{Cosurfaces without underlying manifolds}\label{s:nodimM}
Let us now consider a Hilbert space $H.$ 
Mimicking section \ref{s:cosurface}, we show how
the notions can be extended without the codimension 1 assumption.
\subsection{Settings}
Let $M$ be a $d-$dimensional connected oriented
Riemannian manifold. Let $\mch_{\vee,n}$ (resp. $\mch_{*,n}$ be the set of
embedded, oriented,
smooth, closed, connected  $n$-submanifolds of $H$
with piecewise smooth border (resp. the set of
embedded, oriented,
smooth, closed, connected  $n$-submanifolds of $H$
with piecewise smooth border and with initial and final parts). Let us recall the
assumptions for the $\vee$ composition:

$\bullet$ Let $(s_1,s_2)\in \mch_{\vee,n}\times\mch_{\vee,n}.$
If
\begin{enumerate}
\item $s_1 \cap s_2 \subset \partial s_1 \cap \partial s_2$ is a non empty
$(d-2)$ piecewise smooth manifold and

\item the orientations induced on $s_1 \cap s_2$ by $s_1$ and $s_2$
are opposite.
\end{enumerate}
then we define $s_1 \vee s_2$ by gluing $s_1$ and $s_2$ along their common border.
The orientation of $s_1\vee s_2$
is the one induced by $s_1$ and $s_2,$
and by the gluing conditions, if $s_1 \vee s_2$ is itself a (non piecewise) smooth
submanifold of $H$, then there is an orientation on $s_1\vee s_2$
that generates the orientation of $s_1$ and of $s_2.$

We denote by $\Sigma_{\vee,n}$ the corresponding space 
of piecewise smooth $n$-submanifolds of $H,$
and set $$\Sigma_{\vee,H}= \coprod_{n \in \N} \Sigma_{\vee,n}.$$
Notice that we can extend $\vee$ componentwise to $\Sigma_{\vee,H}.$

\vskip 12pt
\noindent

Let us do the same for $*$.
$\mch_{*,n}$ is the set of (unparametrized, but oriented) $n$-submanifolds $s$ of $H$
 equipped in addition with a prescribed orientation of smooth components of
$\partial s .$
Initial parts of $\partial s$, noted $\alpha(s),$ are those for which the
prescribed orientation is opposite to the one induced by $s,$
and the final parts, noted $\beta(s),$
are the ones for which they co\"incide.
These are exactly the definitions given in section \ref{s:cosurface}, we can define
$\Sigma_{*,n}$ the corresponding space of piecewise smooth submanifolds and we set
$$\Sigma_{*,H}= \coprod_{n \in \N} \Sigma_{*,n}$$ with extension of $*$ componentwise.
In what follows, when it carries no ambiguity,
we omit the notations $\vee$ or $*$ to denote the composition rule by multiplication,
with an adequate choice of Lie group $G$ 
(we choose $G$ to be abelian for $\Sigma_{\vee,H} $).
Anyway, we note the group law of $G$ by multiplication too. 
The following notions can be extended straightaway; 
the only conceptual difference is that the orientation of a surface $s$ 
cannot be compared with the orientation of an underlying manifold $M.$
The manifold $M$ is replaced by the CW-complex obtained by gluing
the domains $A_k$ along the borders.

\begin{Definition}

$\bullet$ A $G$-valued cosurface is a map $$c : \Sigma_H \rightarrow G$$
such that
\begin{enumerate}
\item $\forall (s_1,s_2)\in \Sigma_H \times \h, c(s_1 s_2) =c(s_1)c(s_2)$ and
\item We denote by $\tilde s$ the same $n$-submanifold as
$s\in \Sigma_M$ with opposite orientation on $s$ and $\partial s$. Then
$\forall s \in \Sigma_M, c(\tilde s) = c(s)^{-1}.$
\end{enumerate}
We denote by $\Gamma(G)$ the set of all $G$-valued cosurfaces.

$\bullet$
Let $m \in \mathbb{N}^*.$ An $(m,n)-$\textbf{complex} on $H$ is an $m-$uple $K=(s_1,...s_m) \in \left(\Sigma_{n}\right)^m$ such that $s_i \neq s_j$ for $i\neq j.$ We define
$$ C(K) = (C(s_1),...,C(s_n))$$ where each $C_i \in \Gamma(M;G).$ We denote
by $\mathcal{K}_n$ the set of complexes of dimension $n$ and of any length $m.$

$\bullet$ An $n-$complex $K = (s_1,...,s_n)$ is \textbf{regular} if
  $\forall (i,j) \in {\mathbb{N}_n}^2,$ $$ i \neq j \Rightarrow s_i \cap s_j \subset \partial s_i \cap \partial s_j.$$

$\bullet$ There is a natural notion of subcomplex of a complex $K.$ If $K=(s_1,...,s_n),$ a subcomplex $L$ is a subsequence of $K,$ that is $$\exists l< m\leq n, L = (s_l,...,s_m)=(s_i)_{l\leq i \leq m}.$$

$\bullet$ Let $K$ be a regular $n-$complex. $K$ is called \textbf{saturated}
if and only if there is a $(n+1)-$ complex $A$ such that $\bigcup_{i \in \mathbb{N}_n}s_i$ defines the borders
of the $(n+1)-$surfaces of $A$  by connected and simply connected closed subsets.
In other words, the family $(A_k)_k$ of closed connected and simply connected subsets of $H$ defines, by gluing along $K$, a CW-complex also noted $M$
such that
\begin{enumerate}
\item $\bigcup_{k}A_k = M$
\item for any two indexes $k$ and $k',$
if $k \neq k',$ $$A_k \cap A_{k'} \subset \partial A_k \cap \partial A_{k'} \subset\bigcup_{i \in \mathbb{N}_n}s_i.$$
\end{enumerate}

$\bullet$
We say that a complex $K$ \textbf{splits} through the subcomplex $L = (s_i)_{l \leq i \leq l'}$
if $\bigcup_{s \in L} s $ splits $M$ (as a subset of $H$) into two topologically
connected components $M^+$ and $M^-$ such that $K \cap M^+=(s_i)_{i>m}$ and $K \cap M^- =(s_i)_{i>l'}$
\end{Definition}

\subsection{Dimension extension}

We use here the idea of the previous defintion, gluing together simply connected $(n+1)-$surfaces $A_i$ along a $n-$complex $K$ 
in order to get, by induction on the dimension $N$ of the surfaces, 
a construction of $n-$cosurfaces from lower dimensions.

\begin{Definition}
A complex $K\in \mathcal{K}_n$ is called \textbf{saturated} if it can be embedded
into a $(n+1)$-submanifold $M$ of $H,$ $M \in VCob_{m+1},$ for which
it is a complex for cobordism. $M$ is called a \textbf{saturation} of $K.$
\end{Definition}

In this definition, we then obtain $M$ by gluing along $K$ a family of elements of
$\Sigma_{n+1}$ which are connected and simply connected. 

\begin{Definition}
A complex $K\in \mathcal{K}_{\vee,n}$ is called \textbf{weakly saturated} if there is 
a complex $A=(A_1,...,A_k) \in \mathcal{K}_{*,n+1}$ obtained 

- topologically by gluing each set $A_i$ on $L \subset K$ with respect to the borders $\partial A_i$
such that there is a bijective map $\partial A_i \rightarrow L.$ 

- the orientations of the borders on $\partial A_i$ correspond to the orientations of the submanifolds in the sequence $L.$    

- any $s \in K$ is at least glued once.

$A$ is called a \textbf{weak saturation} of $K.$
\end{Definition}

Notice that with this definition, the orientation of each $A_i$ is left free of choice.
Moreover, the definition can be extended straightaway by 
replacing $\vee$ by $*$ since the initial and final parts of $s \in K$
do not interfer with the gluing.

\begin{Definition} \label{cosurface:dim+1}
Let $c_K$ be defined on the surfaces $A \in \Sigma_{*,n+1}$ that are smooth with 
piecewise smooth border $\partial A \subset K.$
Then we define $$c_K(A) = \prod_{s \in \partial A \subset K} c\circ \varphi_A(s)$$ 
(the product is with respect to the order in $K$).
\end{Definition}

\begin{Proposition} \label{abel}
If  $G$ is abelian, $c_K$ is the restriction of a cosurface $c'$ on $\Sigma_{*,n+1}$ that co\"incide
with $c_K$ on any $A$ where it is defined.   
\end{Proposition}

\vskip 12pt
\noindent
\textbf{Proof.}
First, for a complex $K_1$ such that $K \prec K_1,$ we can get the values 
of $c$ on $K_1.$
Then,    
all we have to show is that, given $A$ an $(n+1)$-submanifold of 
$H$ with piecewise smooth border
along $K$ and taking $A' = (A_1,...A_l)$ such that 
$A \prec A',$ in other words $A = A_1*...*A_l,$
taking $K'$ the (unordered) skeleton of the gluing $A_1*...*A_l,$ one has 
$$c_K(A) = c_{K'}(A_1)...c_{K'}(A_l).$$
There is of course an ambiguity on the order of $K'$ but since 
$G$ is abelian, for each $s \in K'$, 
we have only to count the number of indexes $j\in \N_l$ such that 
$\varphi_{A_j}(s)= s$ and compare it to the number
of indexes such that $\varphi_{A_j}(s) = \tilde{s}.$
Since we are in the $*$-composition, we have at each edge 
$s \not\subset K$ only one index of each type, which are the $A_j$'s for which $s$ 
is in the initial part and the final part respectively.
So that the contributions that are ``interior'' compensate. \qed 

\begin{Definition}
Let $K\in \mathcal{K}_n.$ An \textbf{overcomplex} $K'$ of $K$ is a weakly 
saturated complex such that there exists a complex $K_1$ with $K \prec K_1$ and
 $K_1 \subset K'$ with preserved order. 
\end{Definition}

\begin{Proposition} \label{extnonabelian}
If $G$ is non abelian, $c_K$ is the restriction 
of a cosurface $c'$ on $\Sigma_{*,n+1}$ such that:
for any overcomplex $K'$  of $ K,$ for any $s\in K'\backslash K,$ 
there exists $c''$ a cosurface
on $\Sigma_{*,n+1},$  $$c'(s) = c''_K(s) \in Z(G).$$  
\end{Proposition}  

\vskip 12pt
\noindent
\textbf{Proof.}
Let $A$ be a complex as in definition \ref{cosurface:dim+1}. 
The main points are to know

- what happens for $K_1$ such that $K \prec K_1$? Here again, we can make use of the complex $c$ (we recall that if $K \prec K_1,$ the order is preserved). 

- What happens for an overcomplex K' that the skeleton of a complex $A'$ such that $A$? 
what to do with the "interiors"? 
Let $K_1$ be the complex made of elements $s \in K'$ 
Let is assign arbitrarily a value $c''(s) = g \in Z(G)$ (e.g. take $g=e$)  
if $s\in K'-K_1.$  
We have
\begin{eqnarray*} 
\prod_{i = 1}^l c_K(A_i)& = & \prod_{i = 1}^l\prod_{s \in \partial A_i \subset K} c\circ \varphi_{A_i}(s) \\
& = & \prod_{i=1}^l\prod_{s \in \partial A \subset K_1} c\circ \varphi_{A_i}(s)
\end{eqnarray*}
on one hand. And on the other hand, 
\begin{eqnarray*} 
\prod_{j = 1}^{l'} c_K(A'_j)& = & \prod_{i = 1}^{l'}\prod_{s \in \partial A'_j \subset K} c\circ \varphi_{A_j}(s) \\
& = & \left(\prod_{j=1}^{l'}\prod_{s \in \partial A \subset K_1} c\circ \varphi_{A_i}(s)\right)
\left(\prod_{j=1}^{l'}\prod_{s \in \partial A \subset K-K_1} c\circ \varphi_{A_i}(s)\right)\\
& = & \left(\prod_{j=1}^{l'}\prod_{s \in \partial A \subset K_1} c\circ \varphi_{A_i}(s)\right).e
\end{eqnarray*}
for the same reasons as in the last proof. \qed

\vskip 12pt
\noindent
\textbf{Example: Dimension extension of the holonomy cosurface}   

We have given a way to extend cosurfaces to higher dimensions. 
Let us now use it to show that the notion of non abelian cosurface 
is not void on $\Sigma_{*,n}$ for $n>1.$
For this, let us consider the infinite lattice 
$$\Z^\infty =\left\{ (u_n)_{n \in \N} \in \Z^\N | u_n \neq 0 \hbox{for a finite number of indexes } \right\}$$
Let $G$ be a unimudular Lie group, and let $\theta$ be a connection on $H \times G.$
Along the edges of this lattice in particular and more generally on any path, 
we can define the holonomy cosurface as in section \ref{sss:holcosurface}. 
Let $K$ be a weakly saturated $1-$ complex along the edges of $\Z^\infty.$
Let us now consider the squares that are described by this lattice, 
and more especially those that are gluing along $K.$ 
They are $2-$submanifolds on $H$ and we can define a cosurface $c_{2,K}$ on 
$\Sigma_{2,*}$ that reads as in theorem \ref{extnonabelian}. 
Then, for $K$ ``large enough'', we can find a weakly saturated complex of squares 
$ A = (A_1,...,A_n)$ to reproduce the procedure to get a non abelian 
cosurface on cubes, and this until we reach the dimension $d$ 
for which there no longer exists any weakly 
saturated complex ``based'' on $K.$

For example, take the cube ABCDEFGH and take a complex made of 
its 12 (oriented) segments and $\theta$ a $SU(N)$-connection.
With this choice, $Z(G)=e$ so that $c_{2,K}$ is uniquely determined.      
Then $c_{2,K}$ is non abelian on the 6 faces of the cube, and by 
choosing an order on the faces, i.e. by choosing  complex A made of its 6 faces, 
$c_{2,K}$ extends to a cosurface $c_{3,A,K}$ which is trivial except on 
the cube viewed as an element of $\Sigma_{3,*}.$ If one wants to get 
into another dimension (e.g. on the hypercube), 
the cosurfaces obtained are trivial except for the chosen manifolds along 
the cube.

\section{Cosurfaces and Cobordisms}
\label{s:VCob}

 By `manifold' $M$ we always mean a smooth finite dimensional manifold, 
possibly with boundary; if the boundary, denoted $\partial M$, is void, 
the manifold is said to be {\em closed}. If $N$ is a an oriented manifold 
we denote by $N^{-}$ the manifold $N$ with the opposite orientation.

Let $X_1$ and $X_2$ be oriented  closed submanifolds, both of dimension $d-1$, where $d$ is a positive integer.  By a   \textbf{pre-cobordism} $(Y,\phi_1,\phi_2):X_1\to X_2$ we mean
an oriented manifold $Y$    along with an orientation preserving diffeomorphism
$$\phi:X_1^{-}\sqcup X_2 \to \partial Y.$$
  An isomorphism from a pre-cobordism $Y:X_1\to X_2$ to  a pre-cobordism $(Y',\phi'_1,\phi'_2):X_1\to X_2$ is an orientation-preserving diffeomorphism $f:Y\to Y'$ such that $f\circ\phi_1=\phi'_1$ and $f\circ\phi_2=\phi'_2$.
  A \textbf{cobordism} is a pre-cobordism up to isomorphisms.

Next, if $Y:X_1\to X_2$ is a cobordism, and $Y$ is equipped with a top-dimensional volume form (a measure of volume) we say that $Y$ is a {\em volume} pre-cobordism.
A \textbf{volume cobordism} is then an equivalence class of such 
pre-cobordisms, where the equivalence relation is obtained by using only 
orientation preserving and volume-preserving diffeomorphisms. A theorem of Morse 
guarantees that any two diffeomorphic  compact oriented manifolds of equal volume 
are diffeomorphic by means of a volume-preserving diffeomorphism, and so to 
restrict the considerations to volume cobordism is not a huge restriction).  These 
outlines are  based on a work by Thierry L\'evy \cite{LevyvQFT2011}.

To keep notation under control, we will simply think of a cobordism from $X_1$ to 
$X_2$ as an oriented manifold $Y$, of dimension $d$, 
running ``from'' $X_1$ ``to'' $X_2$. 
Composition of cobordisms is defined in the natural way.   The ``identity''
cobordism $X\to X$  is given by the oriented manifold 
$X\times [0,1]$ along with the mapping
$$X^{-}\sqcup X \to X\times [0,1]$$
which takes $p\in X^{-}$ to $(p,0)$ and $p\in X$ to $(p,1)$.
Let ${\mathbf {VCob}}_d$ be the category whose objects are $d-1$ 
dimensional closed  oriented manifolds and whose morphisms are volume cobordisms.

\vskip 12pt
We may also work within a fixed oriented $d$-dimensional manifold $M$,  
equipped with a volume form, and operate only with cobordisms which are 
(full dimensional) submanifolds of $M$.  Let ${\rm VCob}_M$ be the set of 
all such cobordisms.

\subsection{Adapted saturated complexes}

Consider $Y$ and $Y'$ two morphisms
in the category $VCob_d,$ seen as two $d$-dimensional manifolds equipped with
their borders and volume form.
Assume also that $\alpha(Y)=\beta(Y')$ so that $Y\circ Y'$ exists in $VCob_d.$
Fix now a regular saturated complex $K''$ in $Y\circ Y'$ and for the
$*$-construction so that the set $K''\cap\alpha(Y)$ is made of
complexes on each connected component of $\alpha(Y').$ What we want to construct is
a composition rule for saturated complexes adapted to the composition of
morphisms in $VCob_d.$ Namely, we want to build two complexes
$K\subset Y$ and $K' \subset Y'$ and a ``composition rule'' based on
the composition $*$ for which $K * K' = K''.$

\begin{Definition}
Let $Y\in Mor(VCob_d)$ and let $K$  b a regular saturated complex of $Y.$
Then $K$ is \textbf{adapted} if
for each $ x \in \partial Y \cap K $ and for each$ s \in K \hbox{ such that } x \in s,$
$$x \in \partial s$$ and $$x \in \alpha(s) \Leftrightarrow x \in \alpha(Y).$$
\end{Definition}






Intuitively, the complex $K$ is adapted if it satisfies a property of transversality on the border of $Y,$ and if 
the initial and final prts of $Y$ coincide with the corresponding initial and final parts of the surfaces of $K$ reaching the border.
Now, let $$Y'' = Y \circ Y'.$$
We want to build a regular saturated complex of $Y''$ that splits into $Y$ and $Y'.$ For this, we need 3 parts:  
$$ \left\{ \begin{array}{lcl} K_1 & = & \{\sigma \in K'' | \sigma \subset \alpha(Y)=\beta(Y') \}\\ K_0 & = & \left\{\sigma \in K'' | \sigma \subset Y \hbox{ and } \sigma \notin K_1\right\},\\
K'_0  &= &  \left\{\sigma \in K'' | \sigma \subset Y' \hbox{ and } \sigma \notin K_1 \right\}
 \end{array}\right.$$
such that $K_1$ is a covering of $\alpha(Y),$ and $K_0$ and $K'_0$ are adapted saturated complexes of $Y$ and $Y'$ respectively. Notice that under this condition, $\alpha(Y)$ is a $m-1$-manifold, and the borders $\partial s$, with $s \in K_1$, can define a complex on $\alpha(Y)$ by their smooth components up to re-ordering.
This is what we precise first, and then we give a precise construction from cutting and pasting.

\subsection{Border reduction} \label{ss:borderVCob}
Let $K$ b an adapted complex on $Y \in Mor(VCob_d).$
Let $A$ be the covering of $Y$ with respect to $K$
and let us consider $A_k \in A$ such that
$\partial A_k \cap \alpha(Y) $ has a non empty interior in
$\partial Y$ (one can replace here $\alpha(Y)$ by $\beta(Y)$).
Let $\tilde A_k$ be the closure in $\partial Y$ of the interior
of a connected component of $\partial A_k \cap \partial Y.$
This is a connected subset of $\partial Y, $ not necessarily simply connected.

$\bullet$ \underline{Orientation} $\tilde A_k$ is a $(d-1)$ manifold
with boundary, with
the orientation induced by the orientation of the border of $Y.$

$\bullet$ \underline{Initial and final parts} Now, let us consider $\partial \tilde A_k.$ This is a $(d-2)$ piecewise
smooth manifold, since it is a subset of $\bigcup_{s \in K} \partial s.$
Let $s \in K\cap \partial A_k$ such that
$s \cap \tilde A_k \subset \partial \tilde A_k.$
Then the orientation on $s \cap \partial \tilde A_k$ is the one
induced by the orientation of $s,$ which defines whether
it is an initial or a final part.

Notice that we have here no induced order from the adapted complex $K$ to
the border reduction. The \textbf{border reduction} is a non ordered
regular complex on
$\partial Y,$ which is non necessarily saturated because it defines a partition
of $\partial Y$ into subsets which are non necessarily simply connected, 
with orientations induced by $Y$ and $K.$

\subsection{Complexes for cobordism, cosurfaces and measures part I: cutting}

We now give a more restricted class of complexes.
\begin{Definition}

We say that $K''$ is an $n$-\textbf{complex for cobordism} if $K'' = K_a'' \cup K_\alpha'' \cup K_\beta''$
if \begin{itemize}
\item $K''_\alpha$ is a covering of $\alpha(Y'') = \alpha(Y')$
\item $K''_\beta$ is a covering of $\beta(Y) = \beta(Y'')$
\item $K''_a $ is a saturated complex of $Y'',$
\end{itemize} 
\end{Definition}
We now need to say how we ``cut'' $Y'' \in Mor(VCob_d).$
Let $Y, Y'\in Mor(VCob_d)$ such that $Y'' = Y \circ Y'$ exists. We say that we can cut $(Y'', K'')$ if there exists 
$\theta \in \mathcal{P}(\N_n)$ such that $$K_a'' = \left\{ \sigma_i \in K''| i \in \theta \right\}$$
is an adapted complex in $Y'$ that splits into  $K''_a = K_a \cup K'_a \cup K_b''$, where $K_a$ and $K'_a$ are adapted complexes of $Y$ and $Y'$ respectively, and 
$$K_b'' = \left\{ \sigma_i \in K''| i \not\in \theta \right\}$$
which defines a covering of $\alpha(Y)$

The sets $K_a, K'_a, K''_\alpha, K''_\beta$ and $K_b$ are equipped with the order induced by $K'',$ and gathering the correxponding parts, we get two complexes for cobordism

- $K = K_a \cup K_b \cup K_\beta''$ on $Y$

- $K' = K'_a \cup K_\alpha \cup K_b$ on $Y'.$ 

\begin{rem}
$K_b,$ as a subcomplex of $K'',$ splits $M.$ 
\end{rem}

Now let us turn to measures. For this,
we now take a stochastic cosurface $C''_N$ on $Y''$ adapted to the cobordism,
that is one that can be divided into two stochastic cosurfaces $C_N$ and $C'_N$
on $Y$ and $Y'$ which coincide on
$$\alpha(Y) \cap \Sigma_*=
\left\{ \sigma \in \Sigma_*| \sigma \subset \alpha(Y)  \right\}.$$
Since the order on $K''$ determines the orders on subcomplexes, for each domain $A''_i$
we define $\phi_{A_i}\circ C''$ that equals to
$\phi_{A_i}\circ C$ or $\phi_{A_i}\circ C'$
(we recall that we have $A_i \subset Y$ or $A_i \subset Y'$ since $K''$
is a complex for cobordism)
and each domain is connected and simply connected. If $G$ is non abelian, we assume that the 
indexation of the family $(A_k)_k$ such that the indexes of the domains in $Y'$ are in the beginning of the list, 
and that the indexes of the domains in $Y$ are in the end of the list. If $G$ is abelian, this assumption is not necessary.
\begin{Theorem} \label{muVCob}
$$\mu_K.\mu_{K'} = \mu_{K''}.$$
\end{Theorem}
\vskip 12pt
\noindent
\textbf{Proof.}
Let us build two groups in the formula of definition
\ref{muK}, namely, with the notations of Theorem \ref{muVCob},
$$\mu_{K''}^Q(c'') = \prod_{i =1}^k q_{|A_i|}\left(\phi_{A_i}\circ c''(K)\right).$$
This formula does not depend on the order among the indexes $\N_k,$ so that
we can define a twofold partition ${I,J}$ of $\N_k$ defined as follows: $I$ (resp. $J$)
is the set of indexes $i$ such that $A_i \subset Y$ (resp. $A_i \subset Y'$).
Then
\begin{eqnarray*}
\mu_{K''}^Q(c'') &=& \left(\prod_{i \in I} q_{|A_i|}\left(\phi_{A_i}\circ c''(K)\right)\right).\left(\prod_{j \in J} q_{|A_j|}\left(\phi_{A_j}\circ c''(K)\right)\right)\\
& = & \left(\prod_{i \in I} q_{|A_i|}\left(\phi_{A_i}\circ c(K)\right)\right).\left(\prod_{j \in J} q_{|A_j|}\left(\phi_{A_j}\circ c'(K)\right)\right)\\
& = & \mu_K^Q(c).\mu_{K'}^Q(c') \qed \end{eqnarray*}

\begin{rem}
In order to get saturated complexes we had to add a complex on the border of the manifold $Y.$ This assumption was
not explicitely present in the papers \cite{AHK,AHKH1,AHKH3,AHKH4} where open manifolds were 
also considered. For volume cobordism, 
only compact manifolds with boundary are considered.
The link with finite volume open manifolds can be done with the particular case of
cosurfaces $c$ such that, for any complex for cobordism $K$,
we have the property $c|_{K_b}=e.$
\end{rem}
\subsection{Complexes for cobordism, cosurfaces and measures part II: pasting}
Now let us consider the inverse problem, and let us only point out extra facts that give
``anomalies'' to the pasting procedure.
Let us consider $Y,Y' \in VCob_d,$ equipped with two complexes for cobordisms $K$ and $K'$
and two cosurfaces $c$ and $c'.$ Here are the conditions to
be able to build up a complex for cobordism $K''$on $Y''=Y\circ Y'$:

\vskip 5pt
(A) $K_b\cap \alpha(Y) = K'_b\cap \beta(Y')$ with corresponding orientations, initial and final parts on each hypersurface and on each border.
\vskip 5pt
(B) $c|_{K_b \cap \alpha(Y)} = c'|_{K'_b\cap \beta(Y')}.$

With this, one can build up $c''$, but one cannot build up $K''$ in an unified way.
This depends on a choice of reindexation, 
compatible with the orders of $K_a,$ $K_b,$ $K'b$ and $K'a$
that we have recovered by ``extraction'' from $K''.$

\begin{Proposition}
There exists such a cobordism $K''$ with the properties (A) and (B).
\end{Proposition}

\noindent
\textbf{Proof.}
Let us start with $$K''_0 = K_b \cap \alpha(Y) = K'_b \cap \beta(Y') = (s_1,...,s_l).$$
We build up by induction a complex $K''$ which satisfies (A) and (B).
In a complex $K_e$, for $s \in K_e,$ we note by $bef_{K_e}(s)$ the subcomplex of
elements of $K_e$ before $s$ in the list, and by $aft_{K_e}(s)$ the subcomplex of
elements that are after $s$ in the list.

\vskip 5pt
$\bullet$ \underline{First step}

Let $$K_1'' = bef_K(s_1) \cup bef_{K'}(s_1)\cup$$
(this union is an ordered union, made first of the ordered set $bef_K(s_1)$, secondly of $bef_{K'}(s_1)$ and finally of $K_0$).

\vskip 5pt
$\bullet$ \underline{Intermediate steps}

Let $i \in \N_{l-1}.$ Assume that we know $K_i''.$
Set
\begin{eqnarray*}K_{i+1}''&= &bef_{K_i''}(s_{i+1}) \cup \left(aft_{K}(s_i)\cap bef_{K}(s_{i+1}) \right) \cup \\
&&\left(aft_{K'}(s_i)\cap bef_{K'}(s_{i+1})\right) \cup aft_{K_{i}''}(s_i).\end{eqnarray*}
(with ordered union)

\vskip 5pt
$\bullet$ \underline{Final step}
We have obtained $K''_l$ which last element is $s_l.$
Then $$K'' = K''_l \cup aft_{K}(s_l) \cup aft_{K'}(s_l).$$

\vskip 5pt
Then one can extract $K$ and $K'$ from $K''$ with the desired order. \qed

\vskip 12pt
\noindent
Since the corresponding coverings 
$A=(A_i)_{i \in I}$ of $Y$ and $A'=(A'_j)_{j \in J}$ of $Y'$ are 
well defined and since all the quantities
depend only on the indexation of the hypersurfaces 
on the borders of each domain, with
an order already defined by $K$ and $K'$ and that 
\textit{will not be changed} while passing to $K'',$
the quantity
$$\left(\prod_{i \in I} q_{|A_i|}\left(\phi_{A_i}\circ c(K)\right)\right).\left(\prod_{j \in J} q_{|A_j|}\left(\phi_{A_j}\circ c'(K)\right)\right) = \mu_K^Q(c)\mu_{K'}^Q(c')$$
corresponds to the (classical) definition of $\mu_{K''}^Q(c'')$ for any possible choice of indexation for $K''.$

\section{Algebras and groups of series: applications to cobordism and complexes}
Let us now gather the framework of cobordism and stochastic cosurfaces into series. First, fix $m>1$ the dimension of the theory of cobordism.
The set of indexes is the one described in Theorem \ref{cob}, namely, 
$\Gamma \subset \coprod_{p \in \N^*} I_p$, resp. $\Gamma(N) \subset \coprod_{p \in \N^*} I_p(N)$, be a family of indexes, stable under $*,$ such that $\forall p \in \N^*,$
\begin{enumerate}
\item  $\forall p \in \N^*,$ $\Gamma \cap I_p$ is finite or, more generally;
\item $\forall \gamma \in \Gamma,$ the set of pairs $(\gamma',\gamma'') \in \Gamma^2$ such that $\gamma = \gamma' *`\gamma''$ is finite.
\end{enumerate}
We fix a family $$\A_\Gamma = \prod_{\gamma \in \Gamma}  \A_\gamma $$
of diffeological vector spaces $\A_\gamma$ that are regular (with respect to classical integration of paths)
which is equipped with a multiplication $*$ which is defined componentwise $$\A_\gamma \times \A_{\gamma'} \rightarrow \A_{\gamma * \gamma'}$$
diffeologically smooth, and such that  $(\A_\Gamma, +, *)$ is an algebra.

Let us now turn to our motivating example.
Let $\Gamma$ be a family of piecewise smooth manifolds, made along the infinite lattice $\Z^\infty,$ of fixed dimension $m.$
Let us also normalize the volume of a $m-$cube to 1. 
The family $\Gamma$ is $\N-$graded by the volume, and assume that we have a stochastic cosurface on o $(m-1)-$cubes, either defined directly, or by dimension extension.
Assume now that the complex $K,$ supporting the family $\Gamma$, is a (maybe infinite) complex for cobordism. 
Then, following Theorem \ref{muVCob}, we have the map 
$$ \mu_K : M \in \Gamma \mapsto \mu_K(M) \in \A \subset M(\Omega),$$
where $M(\Omega)$ is the space of measures on $\Omega, $ and $\A$ is a complete vector space of measures such that convolution is associative. 
Then, we are in the context of application of the main theorems of this paper, and we have:
\begin{Theorem}
For a fixed choice of stochastic cosurface, the mapping $\gamma \in \Gamma \mapsto \mu_K(Gamma)$ defines an element of $ \A_\Gamma$
\end{Theorem} 
The proof is obvious.

\end{document}